\documentclass[12pt,reqno]{article}

\usepackage[usenames]{color}
\usepackage{amssymb}
\usepackage{graphicx}
\usepackage{amscd}

\usepackage[colorlinks=true,
linkcolor=webgreen,
filecolor=webbrown,
citecolor=webgreen]{hyperref}

\definecolor{webgreen}{rgb}{0,.5,0}
\definecolor{webbrown}{rgb}{.6,0,0}

\usepackage{color}
\usepackage{fullpage}
\usepackage{float}

\usepackage{graphics,amsmath,amssymb}
\usepackage{amsthm}
\usepackage{amsfonts}
\usepackage{latexsym}
\usepackage{epsf}

\begin{document}

\begin{center}
\vskip 1cm{\LARGE\bf A proof of Catalan's Convolution formula
}
\vskip 1cm
\large
Alon Regev\\
Department of Mathematical Sciences\\
Northern Illinois Univeristy\\
DeKalb, IL\\
\href{mailto:regev@math.niu.edu}{\tt regev@math.niu.edu} \\
\end{center}

\vskip .2 in

\begin{abstract}
We give a new proof of the $k$-fold convolution of the Catalan numbers. This is done by enumerating a certain class of polygonal dissections called $k$-in-$n$ dissections. Furthermore, we give a formula for the average number of cycles in a triangulation.
\end{abstract}

\theoremstyle{plain}
\newtheorem{theorem}{Theorem}
\newtheorem{corollary}[theorem]{Corollary}
\newtheorem{lemma}[theorem]{Lemma}
\newtheorem{proposition}[theorem]{Proposition}

\theoremstyle{definition}
\newtheorem{definition}[theorem]{Definition}
\newtheorem{example}[theorem]{Example}
\newtheorem{conjecture}[theorem]{Conjecture}

\theoremstyle{note}
\newtheorem{note}[theorem]{Note}

\newcommand{\NN}{\mathbb{N}}
\newcommand{\ZZ}{\mathbb{Z}}
\newcommand{\QQ}{\mathbb{Q}}
\newcommand{\im}{\mbox{\upshape Im}}
\newcommand{\summ}{\sum\limits}

\newtheorem{Theorem}{Theorem}[section]
\newtheorem{Proposition}[Theorem]{Proposition}
\newtheorem{Corollary}[Theorem]{Corollary}
\theoremstyle{definition}
\newtheorem{Example}[Theorem]{Example}
\newtheorem{Remark}[Theorem]{Remark}
\newtheorem{Problem}[Theorem]{Problem}

\section{Introduction}

The Catalan numbers are defined as follows.

\begin{definition}
For any $n\ge 0$, $$C_n={1\over n+1}{2n\choose n}.$$
For $n<0$, $C_n=0$.
\end{definition}

The Catalan $k$-fold convolution formula is due to Catalan.

\begin{theorem}\cite{C}\label{catcon}
Let $1 \le k \le n$. Then
\begin{equation}\label{catconeq}
\summ_{i_1+\ldots+i_k=n}C_{i_1-1}\cdots C_{i_k-1}={k\over 2n-k}{2n-k\choose n}.
\end{equation}
\end{theorem}
Catalan's original proof \cite{C, FL, GL, L} uses Lagrange inversion. Gessel and Lacrombe \cite{GL} give two proofs which use hypergeometric identities. Tedford \cite{T} exhibits several interpretations of the left-hand side of \eqref{catconeq}. In this note we use another such interpretation, in terms of dissections of polygons, to give a new proof of Theorem \ref{catcon}. We arrive at this proof using Theorem \ref{k-in-n-formula}, which enumerates a class of polygonal dissections called $k$-in-$n$ dissections. As another consequence of this enumeration, in Corollary \ref{ave} we give a formula for the average number of cycles in a triangulation.

\begin{figure}
\begin{center}
\epsfxsize=2.5in
\leavevmode\epsffile{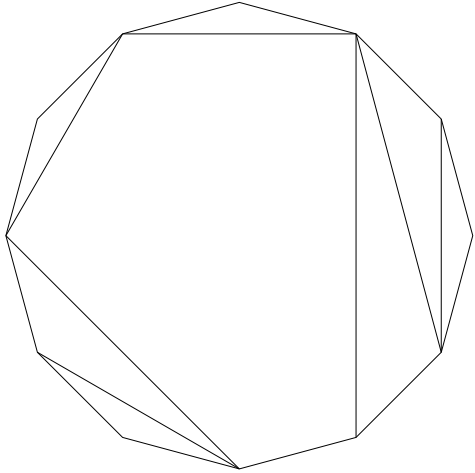}
\end{center}
\caption{Example of a $5$-in-$12$ dissection}
\end{figure}

\section{The $k$-in-$n$ dissections}

\begin{definition}
Let $n \ge 3$ and let $0\le k \le n-3$.
\begin{enumerate}
\item A {\em $k$-dissection} of an $n$-gon is a partition of the $n$-gon into $k+1$ parts by $k$ noncrossing diagonals.
\item A {\em triangulation} of an $n$-gon is an $(n-3)$-dissection.
\item For $k \ge 4$, an {\em $k$-in-$n$ dissection} is an $(n-k)$-dissection of an $n$-gon into one $k$-gon and $n-k+1$ triangles (see Figure 1). A {\em $3$-in-$n$} dissection is a triangulation with one of its $n-3$ triangles marked.
\item Let $f_k(n)$ be the number of $k$-in-$n$ dissections.
\end{enumerate}
\end{definition}

It is well known that for $n\ge 3$ the number of triangulations of an $n$-gon is $C_{n-2}$.
\begin{lemma}\label{k-in-n-rec}
Let $3\le k \le n$. Then
\begin{equation}\label{k-in-n-eq}
(n-k)f_k(n)=n \summ_{i=2}^{n-k+1}C_{i-1} f_k(n-i+1).
\end{equation}
\begin{proof}
The left-hand side of \eqref{k-in-n-eq} is the number of $k$-in-$n$ dissections, with one of the $n-k$ diagonals marked.
These can also be chosen as follows. Choose one vertex $v$ out of the $n$ vertices,
then choose $2\le i\le n-k+1$. Form the diagonal from $v$ to a vertex which is a distance $i$ from $v$ (proceeding, say, counterclockwise along the edges of the $n$-gon). Mark this diagonal. Now choose a triangulation of the resulting $(i+1)$-gon and a $k$-in-$((n-i)+1)$ dissection of the resulting $((n-i)+1)$-gon. Each such choice results in a unique $k$-in-$n$ dissection with one of the diagonals marked.
\end{proof}
\end{lemma}

Lemma \ref{k-in-n-rec} can be used to enumerate the $k$-in-$n$ dissections.
\begin{theorem}\label{k-in-n-formula}
Let $3\le k \le n$. The number of $k$-in-$n$ dissections is
\begin{equation}\label{k-in-n-eq}
f_k(n)={2n-k-1\choose n-1}.
\end{equation}
\end{theorem}

\begin{note}
There is a bijection between $k$-in-$n$ dissections and $k$-crossing partitions of $\{1,\ldots n\}$, as defined in \cite{BMPRSSS}. Thus Theorem \ref{k-in-n-formula} is equivalent to \cite[Theorem 1]{BMPRSSS}.
\end{note}

Theorem \ref{k-in-n-formula} implies the following corollary:

\begin{corollary}\label{ave}
Let $3\le k < n$. The average number of cycles of length $k$ in a triangulated $n$-gon is
$${2n-k-1\choose n-1}{C_{k-2}\over C_{n-2}}.$$
\begin{proof}
Each cycle of length $k$ in a triangulation of an $n$-gon uniquely corresponds to a $k$-in-$n$ dissection together with a triangulation of a $k$-gon. The result then follows from \eqref{k-in-n-eq}.
\end{proof}
\end{corollary}

The following lemmas will be used in the proof of Theorem \ref{k-in-n-formula}. It is well known that for any $n\ge 0$,
\begin{equation}\label{CCeq}
\summ_{i\ge 0} C_{i}C_{n-i}=C_{n+1}.
\end{equation}

\begin{lemma}\label{CC1}
For any $n\ge 1$,
\begin{equation}\label{CC1eq}
\summ_{i\ge 0} i C_{i}C_{n-i}={2n+1\choose n-1}.
\end{equation}
\begin{proof}
Note that
$$\summ_{i\ge 0} i C_{i}C_{n-i}= \summ_{i\ge 0} (n-i) C_{i}C_{n-i}.$$
Therefore by \eqref{CCeq},
$$\summ_{i\ge 0} i C_{i}C_{n-i}={1\over 2} \summ_{i\ge 0}nC_{i}C_{n-i}={n \over 2}C_{n+1}={2n+1\choose n-1}.$$
\end{proof}
\end{lemma}

\begin{lemma}\label{pq}
Let $1\le q\le p\le 2q-1$. Then
\begin{equation}\label{pqeq}
\summ_{i\ge 0} C_i{p-1-2i\choose q-1-i} = {p \choose q}.
\end{equation}
\begin{proof}
We use induction on $q$. If $q=1$ then $p=1$ and both sides of \ref{pqeq} are equal to $1$.
Now suppose $q\ge 2$. If $p=q$ then both sides are equal to $1$.
If $p=2q-1$ then \eqref{pqeq} follows from \eqref{CCeq} and \eqref{CC1eq}, since
\begin{eqnarray*}
\summ_{i\ge 0} C_i{2q-2-2i\choose q-1-i}
&=&\summ_{i\ge 0} C_i (q-i) C_{q-1-i}\\
&=&q\summ_{i\ge 0} C_i C_{q-1-i}-\summ_{i\ge 0} i C_i C_{q-1-i}\\
&=& qC_q-{2q-1\choose q-2}\\
&=&{2q-1 \choose q}.\\
\end{eqnarray*}
Now suppose $q+1 \le p \le 2q-2$. Note that $q-1\le p-1$ and $p-1 \le 2q-2-1=2(q-1)-1$. Therefore by the induction hypothesis, \eqref{pqeq} holds for $p-1$ and $q-1$.
Also $q\le p-1$ and $p-1\le 2q-3<2q-1$, so that \eqref{pqeq} holds for $p-1$ and $q$.  Thus
\begin{eqnarray*}
{p \choose q} &= &{p-1 \choose q-1}+{p-1 \choose q} \\
&=&\summ_{i\ge 0} {1\over i+1} {2i\choose i}{p-2-2i\choose q-2-i}+\summ_{i\ge 0} {1\over i+1} {2i\choose i}{p-2-2i\choose q-1-i}\\
&=&\summ_{i\ge 0} {1\over i+1} {2i\choose i} {p-1-2i\choose q-1-i}. \\
\end{eqnarray*}
\end{proof}
\end{lemma}

\subsection{Proof of Theorem 3}

\begin{proof}
Fix $k\ge 3$ and proceed by induction on $n$. If $n=k$ then both sides are equal to $1$.
Now let $n\ge k+1$. By Lemma \ref{k-in-n-rec} and by the induction hypothesis,
\begin{eqnarray*}
f_k(n)&=&{n\over n-k}\summ_{i=2}^{n-k-1}C_{i-1}f_k(n-i+1)\\
&=&{n\over n-k}\summ_{i=2}^{n-k-1}C_{i-1}{2(n-i+1)-k-1\choose n-i}\\
&=&{n\over n-k}\left(\summ_{i\ge 1}C_{i-1}{2(n-i+1)-k-1\choose n-i}-f_k(n)\right).\\
\end{eqnarray*}
Solving for $f_k(n)$ and applying Lemma \ref{pq}, with $q=n$ and $p=2n-k$,
$$
f_k(n)={n\over 2n-k}\summ_{i\ge 0}C_{i}{2n-k-2i-1\choose n-i-1}={n \over 2n-k}{2n-k\choose n}={2n-k-1\choose n-1}.
$$
\end{proof}

\section{Proof of the Catalan convolution formula}

The next Lemma gives the relation between the number of $k$-in-$n$ dissections and the Catalan convolution.

\begin{lemma}\label{rel}
Let $3\le k < n$. Then
\begin{equation}\label{releq}
kf_k(n)=n \summ_{i_1+\ldots+i_k=n}C_{i_1-1}\cdots C_{i_k-1}.
\end{equation}
\begin{proof}
The left-hand side of \eqref{releq} is the number of $k$-in-$n$ dissections, with one of the vertices of the $k$-gon marked. These can also be chosen as follows. Choose any vertex $v$ of the $n$-gon. For each vertex $v$, choose $i_1,\ldots, i_k$ such that $i_1+\ldots +i_k=n$. This determines the lengths of the sides of a $k$-gon by starting at $v$ and proceeding, say, counterclockwise. For example, in Figure 1, if $v$ is the bottom vertex then the lengths are $1,4,2,2,3$. For each $1\le r \le k$, there is a resulting $(i_r+1)$-gon sharing one edge of the $k$-gon. Each of these $(i_r+1)$-gon can be triangulated in $C_{i_r-1}$ ways, forming a uniquely determined $k$-in-$n$ dissection with one of the of the $k$-gon marked.
\end{proof}
\end{lemma}

The proof of Theorem \ref{catcon} now follows from Lemma \ref{rel}, since
\begin{multline*}
\summ_{i_1+\ldots+i_k=n}C_{i_1-1}\cdots C_{i_k-1}={k\over  n}f_k(n)
={k\over n}{2n-k-1\choose n-1}
={k\over 2n-k}{2n-k \choose n}.\\
\end{multline*}

\end{document}